\documentclass{statsoc}
\usepackage{amssymb}
\usepackage{amsmath}
\usepackage[english]{babel}
\newtheorem{theorem}{Theorem}[section]
\newtheorem{lemma}[theorem]{Lemma}

\newcommand{\R}{\ensuremath{\mathbb{R}}}
\title[]{Maximum likelihood estimation in a partially observed stratified regression model with censored data}
\author[A. Detais, J.-F. Dupuy]{Detais Am\'elie}
\address{Laboratoire de Statistique et Probabilit\'es, UMR 5219, Institut de Math\'ematiques de Toulouse, Universit\'e Toulouse 3, France.\\ Email: Amelie.Detais@math.ups-tlse.fr}
\author[A. Detais, J.-F. Dupuy]{Dupuy Jean-Fran\c cois\footnote{Corresponding author}}
\address{Laboratoire de Statistique et Probabilit\'es, UMR 5219, Institut de Math\'ematiques de Toulouse, Universit\'e Toulouse 3, France.\\ Email: Jean-Francois.Dupuy@math.ups-tlse.fr}
\begin{document}
\begin{abstract}
The stratified proportional intensity model generalizes Cox's proportional intensity model by allowing different groups of the population under study to have distinct baseline intensity functions. In this article, we consider the problem of estimation in this model when the variable indicating the stratum is unobserved for some individuals in the studied sample. In this setting, we construct  nonparametric maximum likelihood estimators for the parameters of the stratified model and we establish their consistency and asymptotic normality. Consistent estimators for the limiting variances are also obtained.
\end{abstract}

\section{Introduction}
This paper considers the problem of estimation in the stratified proportional intensity regression model for survival data, when the stratum information is missing for some sample individuals.

The stratified proportional intensity model (see \cite{and93} or \cite{mart06} for example) generalizes the usual \cite{cox72} proportional intensity regression model for survival data, by allowing different groups -the strata- of the population under study to have distinct baseline intensity functions. More precisely, in the stratified model, the strata divide the sample individuals into $K$ disjoint groups, each having a distinct baseline intensity function $\lambda_k$ but a common value for the regression parameter.

The intensity function for the failure time $T^0$ of an individual in stratum $k$ thus takes the form
\begin{equation}\label{model}
\lambda_k(t)\exp(\beta'X),
\end{equation}
where $X$ is a $p$-vector of covariates, $\beta$ is a $p$-vector of unknown regression parameters of interest, and $\{\lambda_k(t):t\geq 0, k=1,\ldots,K\}$ are $K$ unknown baseline intensity functions.

A consistent and asymptotically normal estimator of $\beta$ can be obtained by maximizing the partial likelihood function \citep{cox75}. The partial likelihood for the stratified model (\ref{model}) is the product over strata of the within-stratum partial likelihoods (we refer to \cite{and93} for a detailed treatment of maximum partial likelihood estimation in model (\ref{model})). In some applications, it can also be desirable to estimate the cumulative baseline intensity functions $\Lambda_k= \int \lambda_k$. The so-called \cite{bres72} estimators are commonly used for that purpose (see chapter 7 of \cite{and93} for further details on the Breslow estimator and its asymptotic properties).

One major motivation for using the stratified model is that it allows to accomodate in the analysis a predictive categorical covariate whose effect on the intensity is not proportional. To this end, the individuals under study are stratified with respect to the categories of this covariate. In many applications however, this covariate may be missing for some sample individuals (for example, histological stage determination may require biopsy and due to expensiveness, may not be performed on all the study subjects). In this case, the usual statistical inference for model (\ref{model}), based on the product of within-stratum partial likelihoods, can not be directly applied.

In this work, we consider the problem of estimating $\beta$ and the $\Lambda_k, k=1,\ldots,K$ in model (\ref{model}), when the covariate defining the stratum is missing for some (but not all) individuals.  Equivalently said, we consider the problem of estimating model (\ref{model}) when the stratum information is only partially available.

The problem of estimation in the (unstratified) Cox regression model $\lambda(t)\exp(\beta'X)$ with missing covariate $X$ has been the subject of intense research over the past decade: see for example \cite{lin93}, \cite{paiklida}, \cite{paik97}, \cite{chen99}, \cite{mart99}, \cite{pons02}, and the references therein. But to the best of our knowledge and despite its practical relevance, the problem of statistical inference in model (\ref{model}) with partially available stratum information has not been yet extensively investigated. Recently, \cite{duplec} studied the asymptotic properties of a regression calibration estimator of $\beta$ in this setting (regression calibration is a general method for handling missing data in regression models, see \cite{crs} for example). The authors proved that this estimator is asymptotically biased, although nevertheless asymptotically normal. No estimators of the cumulative baseline intensity functions were provided.

In this work, we aim at providing an estimator of $\beta$ that is both consistent and asymptotically normal. Moreover, although the cumulative intensity functions $\Lambda_k$ are usually not the primary parameters of interest, we also aim at providing consistent and asymptotically normal estimators of the values $\Lambda_k(t)$, $k=1,\ldots,K$.

The regression calibration inferential procedure investigated by \cite{duplec} is essentially based on a modified version of the partial likelihood for model (\ref{model}). In this paper, we propose an alternative method which may be viewed as a fully maximum likelihood approach. Besides assuming that the failure intensity function for an individual in stratum $k$ is given by model (\ref{model}), we assume that the probability of being in stratum $k$ conditionally on a set of observed covariates $W$ (which may include some components of $X$) is of the logistic form, depending on some unknown finite-dimensional parameter $\gamma$. 

A full likelihood for the collected parameter $\theta=(\beta, \gamma, \Lambda_{k}; k=1,\ldots,K)$ is constructed from a sample of incompletely observed data. Based on this, we propose to estimate the finite and infinite-dimensional components of $\theta$ by using the nonparametric maximum likelihood (NPML) estimation method. We then provide asymptotic results for these estimators, including consistency, asymptotic normality, semiparametric efficiency of the NPML estimator of $\beta$, and consistent variance estimation.

Our proofs use some techniques developed by \cite{murphy94,murphy95} and \cite{parner} to establish the asymptotic theory for the frailty model.

The paper is organized as follows. In Section \ref{dataassump}, we describe in greater detail the data structure and the model assumptions. In Section \ref{npmldesc}, we describe the NPML estimation method for our setting and we establish existence of the NPML estimator of $\theta$. Section \ref{asymprp} establishes the consistency and asymptotic normality of the proposed estimator. Consistent variance estimators are also obtained for both the finite-dimensional parameter estimators and the nonparametric cumulative baseline intensity estimators. We give some concluding remarks in Section \ref{disc}. Proofs are given in Appendix.

\section{Data structure and model assumptions}\label{dataassump}

We describe the notations and model assumptions that will be used throughout the paper.

All the random variables are defined on a probability space $(\Omega,\mathcal C, \mathbb P)$. Let $T^0$ be a random failure time whose distribution depends on a vector of covariates $X\in \R^p$ and on a stratum indicator $S\in\mathcal K=\{1,\ldots,K\}$. We assume that conditionally on $X$ and $S=k$ $(k\in\mathcal K)$, the intensity function of $T^0$ is given by model (\ref{model}). We suppose that $T^0$ may be right-censored by a positive random variable $C$ and that the analysis is restricted to the time interval $[0,\tau]$, where $\tau<\infty$ denotes the end of the study. Thus we actually observe the potentially censored duration $T=\min\{T^0,\min(C,\tau)\}$ and a censoring indicator $\Delta=1\{T^0\leq \min(C,\tau)\}$. If $t\in[0,\tau]$, we denote by $N(t)=1\{T\leq t\}\Delta$ and $Y(t)=1\{T\geq t\}$ the failure counting and at-risk processes respectively.

Let $W\in\R^m$ be a vector of surrogate covariates for $S$ ($W$ and $X$ may share some common components). That is, $W$ brings a partial information about $S$ when $S$ is missing, and it adds no information when $S$ is observed so that the distribution of $T^0$ conditionally on $X, S$, and $W$ does not involve the components of $W$ that are not in $X$. We assume that the conditional probability that an individual belongs to the $k$-th stratum given his covariate vector $W$ follows a multinomial logistic model:
\begin{eqnarray*}
\mathbb P(S=k|W)=\frac{\exp(\gamma_k'W)}{\sum_{j=1}^K\exp(\gamma_j'W)},
\end{eqnarray*}
where $\gamma_k\in\mathbb R^m$ $(k\in\mathcal K)$. Finally, we let $R$ denote the indicator variable which is 1 if $S$ is observed and 0 otherwise. Then, the data consist of $n$ i.i.d. replicates
\begin{eqnarray*}
\mathcal O_i=(T_i, \Delta_i, X_i, W_i, R_i, R_iS_i), \quad i=1,\ldots,n,
\end{eqnarray*}
of $\mathcal O=(T, \Delta, X, W, R, RS)$. The data available for the $i$-th individual are therefore $(T_i, \Delta_i, X_i, W_i, S_i)$ if $R_i=1$ and $(T_i, \Delta_i, X_i, W_i)$ if $R_i=0$.

In the sequel, we set $\gamma_K=0$ for model identifiability purpose and we note $\gamma=(\gamma_1',\ldots,\gamma_{K-1}')'\in(\R^m)^{K-1}\equiv\R^q$. We also note $\pi_{k,\gamma}(W)=\mathbb P(S=k|W)$, $k\in\mathcal K$. Now, let $\theta=(\beta, \gamma, \Lambda_{k}; k\in\mathcal K)$ be the collected parameter and $\theta_0=(\beta_0, \gamma_0, \Lambda_{k,0}; k\in\mathcal K)$ denote the true parameter value. Under the true value $\theta_0$, the expectation of random variables will be denoted $P_{\theta_0}$. $\mathbb P_n$ will denote the empirical probability measure. In the sequel, the stochastic convergences will be in terms of outer measure.

We now make the following additional assumptions:
\begin{description}
\item[(a)] The censoring time $C$ is independent of $T^0$ given $(S, X, W)$, of $S$ given $(X, W)$, and is non-informative. With probability 1, $\mathbb P(C\geq T^0\geq \tau|S, X, W)>c_0$ for some positive constant $c_0$.
\item[(b)] The parameter values $\beta_0$ and $\gamma_0$ lie in the interior of known compact sets $\mathcal B\subset \R^p$ and $\mathcal G\subset \R^q$ respectively. For every $k\in\mathcal K$, the cumulative baseline intensity function $\Lambda_{k,0}$ is a strictly increasing function on $[0,\tau]$ with $\Lambda_{k,0}(0)=0$ and $\Lambda_{k,0}(\tau)<\infty$. Moreover, for every $k\in\mathcal K$, $\Lambda_{k,0}$ is continuously differentiable in $[0,\tau]$, with $\lambda_{k,0}(t)=\partial\Lambda_{k,0}(t)\slash\partial t$. Let $\mathcal A$ denote the set of functions satisfying these properties.
\item[(c)] The covariate vectors $X$ and $W$ are bounded (\textit{i.e.} $\|X\|<c_1$ and $\|W\|<c_1$, for some finite positive constant $c_1$, where $\|\cdot\|$ denotes the Euclidean norm). Moreover, the covariance matrices of $X$ and $W$ are positive definite. Let $c_2=\min_{\beta\in\mathcal B, \|X\|<c_1}e^{\beta'X}$ and $c_3=\max_{\beta\in\mathcal B, \|X\|<c_1}e^{\beta'X}$.
\item[(d)] There is a constant $c_4>0$ such that for every $k\in\mathcal K$, $P_{\theta_0}[1\{S=k\}Y(\tau)R]>c_4$, and the sample size $n$ is large enough to ensure that $\sum_{i=1}^n 1\{S_i=k\}Y_i(\tau)R_i>0$ for every $k\in\mathcal K$.
\item[(e)] With probability 1, there exists a positive constant $c_5$ such that for every $k\in\mathcal K$, $P_{\theta_0}[\Delta R 1\{S=k\}|T, X, W]>c_5$.
\item[(f)] $R$ is independent of $S$ given $W$, of $(T, \Delta)$ given $(X,S)$. The distribution of $S$ conditionally on $X$ and $W$ does not involve the components of $X$ that are not in $W$. The distributions of $R$ and of the covariate vectors $X$ and $W$ do not depend on the parameter $\theta$.
\end{description}
\textsc{Remark 1.} Conditions (b), (c), (d), and (e) are used for identifiability of the parameters and consistency of the proposed estimators. Condition (d) essentially requires that for every stratum $k$, some subjects are known to belong to $k$ and are still at risk when the study ends. The first assumption in condition (f) states that $S$ is missing at random, which is a fairly general missing data situation (we refer to chapters 6 and 7 in \cite{tsi06} for a recent exposition of missing data mechanisms).
\\\\
\textsc{Remark 2.} We are now in position to describe our proposed approach to the problem of estimation in model (\ref{model}) from a sample of incomplete data $\mathcal O_i$, $i=1,\ldots,n$.
\\
Let $\mathcal S$ denote the set of subjects with unknown stratum in this sample. The regression calibration method investigated by \cite{duplec} essentially allocates every subject of $\mathcal S$ to each of the strata, and estimates $\beta_0$ by maximizing a modified version of the partial likelihood for the stratified model, where the contribution of any individual $i$ in $\mathcal S$ to the within-$k$-th-stratum partial likelihood is weighted by an estimate of $\pi_{k,\gamma}(W_i)$ (for every $k\in\mathcal K$). The asymptotic bias of the resulting estimator arises from the failure of this method to fully exploit the information carried by $(T_i, \Delta_i, X_i, W_i)$ on the unobserved stratum indicator $S_i$.
\\
Therefore in this paper, we rather suggest to weight each subject $i$ in $\mathcal S$ by an estimate of the conditional probability that subject $i$ belongs to the $k$-th stratum given the whole observed data $(T_i, \Delta_i, X_i, W_i)$. This suggestion raises two main problems, as is described below.
\\\\
\textsc{Remark 3.} First, we should note that the suggested alternative weights depend on the unknown baseline intensity functions. Therefore, the modified partial likelihood approach considered by \cite{duplec} can not be used to derive an estimator for $\beta_0$. Next, the statistics to be involved in the score function for $\beta$ will depend on the conditional weights and thus, this score will not be expressible as a stochastic integral of some predictable process, as is often the case in models for failure time data. This, in turn, will prevent us from using the counting process martingale theory usually associated with the theoretical developments in failure time models.
\\\\
To overcome the first problem, we define our estimators from a full likelihood for the whole parameter, that is, for both the finite-dimensional -$\beta$ (and $\gamma$)- and infinite-dimensional -$\Lambda_k$, $k\in\mathcal K$- components of $\theta$. Empirical process theory \citep{vdvw96} is used to establish asymptotics for the proposed estimators.

\section{Maximum likelihood estimation}\label{npmldesc}

In the sequel, we assume that there are no ties among the observed death times (this hypothesis is made to simplify notations, but the results below can be adapted to accomodate ties). The likelihood function for observed data  $\mathcal O_i$, $i=1,\ldots,n$ is given by
\begin{eqnarray}\label{odl}
&&L_n(\theta)=\prod_{i=1}^n \left[\prod_{k=1}^K\left\{\lambda_k(T_i)^{\Delta_i}\exp\left(\Delta_i\beta'X_i-e^{\beta'X_i} \Lambda_k(T_i)\right)\pi_{k,\gamma}(W_i)\right\}^{1\{S_i=k\}}\right]^{R_i}\nonumber\\
&&\hspace{2.5cm}\times\left[\sum_{k=1}^K\lambda_k(T_i)^{\Delta_i}\exp\left(\Delta_i\beta'X_i-e^{\beta'X_i} \Lambda_k(T_i)\right)\pi_{k,\gamma}(W_i)\right]^{1-R_i}.
\end{eqnarray}
It would seem natural to derive a maximum likelihood estimator of $\theta_0$ by maximizing the likelihood (\ref{odl}). However, the maximum of this function over the parameter space $\Theta=\mathcal B \times \mathcal G\times \mathcal A^{\otimes K}$ does not exist. To see this, consider functions $\Lambda_k$ with fixed values at the $T_i$, and let $(\partial\Lambda_k(t)\slash\partial t)|_{t=T_i}=\lambda_k(T_i)$ go to infinity for some $T_i$ with $\Delta_iR_i1\{S_i=k\}=1$ or $\Delta_i(1-R_i)=1$.

To overcome this problem, we introduce a modified maximization space for (\ref{odl}), by relaxing each $\Lambda_k(\cdot)$ to be an increasing right-continuous step-function on $[0,\tau]$, with jumps at the $T_i$'s such that $\Delta_iR_i1\{S_i=k\}=1$ or $\Delta_i(1-R_i)=1$. Estimators of $(\beta_0, \gamma_0, \Lambda_{k,0}; k\in\mathcal K)$ will thus be derived by maximizing a modified version of (\ref{odl}), obtained by replacing $\lambda_k(T_i)$ in (\ref{odl}) with the jump size $\Lambda_k\{T_i\}$ of $\Lambda_k$ at $T_i$.

If they exist, these estimators will be referred to as nonparametric maximum likelihood estimators - NPMLEs - (we refer to \cite{zeng07} for a review of the general principle of NPML estimation, with application to various semiparametric regression models for censored data. See also the numerous references therein). In our setting, existence of such estimators is ensured by the following theorem (proof is given in Appendix):
\begin{theorem}\label{exist}
Under conditions (a)-(f), the NPMLE $\widehat\theta_n=(\widehat\beta_n, \widehat\gamma_n, \widehat\Lambda_{k,n}; k\in\mathcal K)$ of $\theta_0$ exists and is achieved.
\end{theorem}
The problem of maximizing $L_n$ over the approximating space described above reduces to a finite dimensional problem, and the expectation-maximization (EM) algorithm \citep{dlr} can be used to calculate the NPMLEs. For $1\leq i\leq n$ and $k\in\mathcal K$, let $w_i(k,\theta)$ be the conditional probability that the $i$-th individual belongs to the $k$-th stratum given $(T_i, \Delta_i, X_i, W_i)$ and the parameter value $\theta$, and let $Q(\mathcal O_i,k,\theta)$ denote the conditional expectation of $1\{S_i=k\}$ given $\mathcal O_i$ and the parameter value $\theta$. Then $Q(\mathcal O_i,k,\theta)$ has the form
\begin{eqnarray*}
Q(\mathcal O_i,k,\theta)=R_i1\{S_i=k\}+(1-R_i)w_i(k,\theta).
\end{eqnarray*}
In the M-step of the EM-algorithm, we solve the complete-data score equation conditional on the observed data. In particular, the following expression for the NPMLE of $\Lambda_k(\cdot)$ can be obtained by: (a) taking the derivative with respect to the jump sizes of $\Lambda_k(\cdot)$, of the conditional expectation of the complete-data log-likelihood given the observed data and the NPML estimator, (b) setting this derivative equal to 0:
\begin{lemma}
The NPMLE $\widehat\theta_n$ satisfies the following equation for every $k\in\mathcal K$:
\begin{eqnarray*}\label{npmleL}
\widehat\Lambda_{k,n}(t)=\int_0^t\sum_{i=1}^n \frac{Q(\mathcal O_i,k,\widehat\theta_n)}{\sum_{j=1}^n Q(\mathcal O_j,k,\widehat\theta_n)\exp(\widehat\beta_n'X_j)Y_j(s)}\,dN_i(s),\quad 0\leq t\leq \tau.
\end{eqnarray*}
\end{lemma}
The details of the calculations are omitted (note how the suggested weights $w_i(k,\theta)$ naturally arise from the M-step of the EM algorithm). We refer the interested reader to \cite{zengcai05} and \cite{sugi06}, who recently described EM algorithms for computing NPMLEs in various other semiparametric models with censored data.

In the sequel, we shall denote the conditional expectation of the complete-data log-likelihood given the observed data and the NPML estimator by $E_n[\widetilde l_n(\theta)]$.

\section{Asymptotic properties}\label{asymprp}

This section states the asymptotic properties of the proposed estimators. We first obtain the following theorem, which states the strong consistency of the proposed NPMLE. The proof is given in Appendix.

\begin{theorem}\label{csy}
Under conditions (a)-(f), $\|\widehat\beta_n-\beta_0\|$, $\|\widehat\gamma_n-\gamma_0\|$, and $\sup_{t\in[0,\tau]}|\widehat\Lambda_{k,n}(t)-\Lambda_{k,0}(t)|$ (for every $k\in\mathcal K$) converge to 0 almost surely as $n$ tends to infinity.
\end{theorem}

To derive the asymptotic normality of the proposed estimators, we adapt the function analytic approach developed by \cite{murphy95} for the frailty model (see also \cite{chang05}, \cite{korui07}, and \cite{lu08}, for recent examples of this approach in various other models).

Instead of calculating score equations by differentiating $E_n[\widetilde l_n(\theta)]$ with respect to $\beta$, $\gamma$, and the jump sizes of $\Lambda_k(\cdot)$, we consider one-dimensional submodels $\widehat\theta_{n,\eta}$ passing through $\widehat\theta_n$ and we differentiate with respect to $\eta$. Precisely, we consider submodels of the form
\begin{eqnarray*}
\eta\mapsto \widehat\theta_{n,\eta}=\left(\widehat\beta_n+\eta h_\beta,\widehat\gamma_n+\eta h_\gamma, \int_0^\cdot \{1+\eta h_{\Lambda_k}(s)\}d\widehat\Lambda_{k,n}(s); k\in\mathcal K\right),
\end{eqnarray*}
where $h_\beta$ and $h_\gamma=(h_{\gamma_1}',\ldots,h_{\gamma_{K-1}}')'$ are $p$- and $q$-dimensional vectors respectively ($h_{\gamma_j}\in\mathbb R^m$, $j=1,\ldots,K-1$), and the $h_{\Lambda_k}$ $(k\in\mathcal K)$ are functions on $[0,\tau]$. Let $h=(h_\beta, h_\gamma, h_{\Lambda_k}; k\in\mathcal K)$. To obtain the score equations, we differentiate $E_n[\widetilde l_n(\widehat\theta_{n,\eta})]$ with respect to $\eta$ and we evaluate at $\eta=0$. $\widehat\theta_n$ maximizes $E_n[\widetilde l_n(\theta)]$ and therefore satisfies $(\partial E_n[\widetilde l_n(\widehat\theta_{n,\eta})]\slash\partial\eta)\left.\right|_{\eta=0}=0$ for every $h$, which leads to the score equation $S_n(\widehat\theta_n)(h)=0$ where $S_n(\widehat\theta_n)(h)$ takes the form
\begin{eqnarray}\label{score}
S_n(\widehat\theta_n)(h)=\mathbb P_n\left[h_\beta' S_\beta(\widehat\theta_n)+ h_\gamma' S_\gamma(\widehat\theta_n)+\sum_{k=1}^K S_{\Lambda_k}(\widehat\theta_n)(h_{\Lambda_k})\right],
\end{eqnarray}
where
\begin{eqnarray*}
&&S_\beta(\theta)=\Delta X-\sum_{k=1}^K Q(\mathcal O,k,\theta)X\exp(\beta'X)\Lambda_k(T),\\
&&S_\gamma(\theta)=(S_{\gamma_1}(\theta)',\ldots,S_{\gamma_{K-1}}(\theta)')' \mbox{ with } S_{\gamma_k}(\theta)=W\left[Q(\mathcal O,k,\theta)-\pi_{k,\gamma}(W)\right],\\
&&S_{\Lambda_k}(\theta)(h_{\Lambda_k})=Q(\mathcal O,k,\theta)\left[h_{\Lambda_k}(T)\Delta-\exp(\beta'X)\int_0^{T}h_{\Lambda_k}(s)\,d\Lambda_k(s)\right].
\end{eqnarray*}
We take the space of elements $h=(h_\beta, h_\gamma, h_{\Lambda_k}; k\in\mathcal K)$ to be
\begin{eqnarray*}
&&H=\{(h_\beta, h_\gamma, h_{\Lambda_k}; k\in\mathcal K): h_\beta\in\R^p, \|h_\beta\|<\infty; h_\gamma\in\R^q, \|h_\gamma\|<\infty;  \\
&&\hspace{2.5cm}h_{\Lambda_k} \mbox{ is a function defined on } [0,\tau], \|h_{\Lambda_k}\|_v<\infty, k\in\mathcal K\},
\end{eqnarray*}
where $\|h_{\Lambda_k}\|_v$ denotes the total variation of $h_{\Lambda_k}$ on $[0,\tau]$. We further take the functions $h_{\Lambda_k}$ to be continuous from the right at 0.

Define $\theta(h)=h_\beta'\beta+h_\gamma'\gamma+\sum_{k=1}^K\int_0^\tau h_{\Lambda_k}(s)\,d\Lambda_k(s)$, where $h\in H$. From this, the parameter $\theta$ can be considered as a linear functional on $H$, and the parameter space $\Theta$ can be viewed as a subset of the space $l^\infty(H)$ of bounded real-valued functions on $H$, which we provide with the uniform norm. Moreover, the score operator $S_n$ appears to be a random map from $\Theta$ to the space $l^\infty(H)$. Note that appropriate choices of $h$ allow to extract all components of the original parameter $\theta$. For example, letting $h_\gamma=0$, $h_{\Lambda_k}(\cdot)=0$ for every $k\in\mathcal K$, and $h_\beta$ be the $p$-dimensional vector with a one at the $i$-th location and zeros elsewhere yields the $i$-th component of $\beta$. Letting $h_\beta=0$, $h_\gamma=0$, $h_{\Lambda_k}(\cdot)=0$ for every $k\in\mathcal K$ except $h_{\Lambda_j}(s)=1\{s\leq t\}$ (for some $t\in (0,\tau)$) yields $\Lambda_j(t)$.

We need some further notations to state the asymptotic normality of the NPMLE of $\beta_0$. Let us first define the linear operator $\sigma=(\sigma_\beta,\sigma_\gamma,\sigma_{\Lambda_k};k\in\mathcal K):H\rightarrow H$ by
\begin{eqnarray*}
&&\sigma_\beta(h)=P_{\theta_0}\left[2X\Delta \psi(\mathcal O,\theta_0)\sum_{k=1}^KQ(\mathcal O,k,\theta_0) h_{\Lambda_k}(T) \right]\\
&&\hspace{2cm}+P_{\theta_0}\left[\psi(\mathcal O,\theta_0)X\left\{\psi(\mathcal O,\theta_0)X' h_\beta+S_\gamma(\theta_0)'h_\gamma\right\}\right],\\
&&\sigma_\gamma(h)=P_{\theta_0}\left[2S_\gamma(\theta_0)\Delta\sum_{k=1}^KQ(\mathcal O,k,\theta_0) h_{\Lambda_k}(T) \right]+P_{\theta_0}\left[S_\gamma(\theta_0)S_\gamma(\theta_0)'\right]h_\gamma\\
&&\hspace{3cm}+P_{\theta_0}\left[\psi(\mathcal O,\theta_0)S_\gamma(\theta_0)X'\right]h_\beta,\\
&&\sigma_{\Lambda_k}(h)(u)=h_{\Lambda_k}(u)P_{\theta_0}\left[Q(\mathcal O,k,\theta_0) \phi(u,\mathcal O, k, \theta_0) \right] \\
&&\hspace{2cm}+P_{\theta_0}\left[2\phi(u,\mathcal O,k,\theta_0)\sum_{j>k}Q(\mathcal O,j,\theta_0) \left\{h_{\Lambda_j}(u)-e^{\beta_0'X}\int_0^u h_{\Lambda_j}\,d\Lambda_{j,0}\right.\right.\\
&&\hspace{3cm}\left.\left.-\Delta h_{\Lambda_j}(T)+e^{\beta_0'X}\int_0^T h_{\Lambda_j}\,d\Lambda_{j,0}\right\}\right]\\
&&\hspace{2cm}-h_\beta'P_{\theta_0}\left[2X\psi(\mathcal O,\theta_0)Q(\mathcal O,k,\theta_0)e^{\beta_0'X}Y(u)\right]\\
&&\hspace{2cm}-h_\gamma'P_{\theta_0}\left[2S_\gamma(\theta_0)Q(\mathcal O,k,\theta_0)e^{\beta_0'X}Y(u)\right],
\end{eqnarray*}
where $\phi(u,\mathcal O,k,\theta_0)=Y(u)Q(\mathcal O,k,\theta_0)e^{\beta_0'X}$ and $\psi(\mathcal O,\theta_0)=\Delta -\sum_{k=1}^K Q(\mathcal O,k,\theta_0)e^{\beta_0'X}\Lambda_{k,0}(T)$. This operator is continuously invertible (Lemma \ref{infoope} in Appendix). We shall denote its inverse by $\sigma^{-1}=(\sigma_\beta^{-1},\sigma_\gamma^{-1},\sigma_{\Lambda_k}^{-1};k\in\mathcal K)$.

Next, for every $r\in\mathbb N\backslash \{0\}$, the $r$-dimensional column vector having all its components equal to 0 will be noted by $0_r$ (or by $0$ when no confusion may occur). Let $h=(h_\beta, h_\gamma, h_{\Lambda_k}; k\in\mathcal K)\in H$. If $h_\gamma=0$ and $h_{\Lambda_k}$ is identically equal to 0 for every $k\in\mathcal K$, we note $h=(h_\beta, 0, 0; k\in\mathcal K)$. Let $\widetilde \sigma_\beta^{-1}:\R^p\rightarrow \R^p$ be the linear map defined by $\widetilde \sigma_\beta^{-1}(u)=\sigma_\beta^{-1}((u,0,0;k\in\mathcal K))$, for $u\in\R^p$. Let $\{e_1,\ldots,e_p\}$ be the canonical basis of $\R^p$.

Then the following result holds, its proof is given in Appendix.

\begin{theorem}\label{asympn}
Under conditions (a)-(f), $\sqrt n(\widehat\beta_n-\beta_0)$ has an asymptotic normal distribution $N(0,\Sigma_\beta)$, where
\begin{eqnarray*}
\Sigma_\beta=(\widetilde \sigma_\beta^{-1}(e_1),\ldots,\widetilde \sigma_\beta^{-1}(e_p))
\end{eqnarray*}
is the efficient variance in estimating $\beta_0$.
\end{theorem}

\textsc{Remark 4.} Although $\gamma_0$ and the cumulative baseline intensity functions $\Lambda_{k,0}$ $(k\in\mathcal K)$ are not the primary parameters of interest, we may also state an asymptotic normality result for their NMPLEs. This requires some further notations.
\\
\\
Define $\widetilde \sigma_\gamma^{-1}:\R^q\rightarrow \R^q$ by $\widetilde \sigma_\gamma^{-1}(u)= \sigma_\gamma^{-1}((0, u, 0; k \in \mathcal K))$, let $\{f_1,\ldots,f_q\}$ be the canonical basis of $\R^q$, and define $\Sigma_\gamma=(\widetilde \sigma_\gamma^{-1}(f_1),\ldots,\widetilde \sigma_\gamma^{-1}(f_q))$. Finally, let $h_{(j,t)}=(h_\beta, h_\gamma, h_{\Lambda_k}; k\in\mathcal K)$ be such that $h_\beta=0$, $h_\gamma=0$, $h_{\Lambda_j}(\cdot)=1\{\cdot\leq t\}$ for some $t\in(0,\tau)$ and $j\in\mathcal K$, and $h_{\Lambda_k}=0$ for every $k\in\mathcal K, k\neq j$. Then the following holds (a brief sketch of the proof is given in Appendix):

\begin{theorem}\label{asympnpn}
Assume that conditions (a)-(f) hold. Then $\sqrt n(\widehat\gamma_n-\gamma_0)$ has an asymptotic normal distribution $N(0,\Sigma_\gamma)$. Furthermore, for any $t\in(0,\tau)$ and $j\in\mathcal K$, $\sqrt n(\widehat\Lambda_{j,n}(t)-\Lambda_{j,0}(t))$ is asymptotically distributed as a $N(0, v_j^2(t))$, where
\begin{eqnarray*}
v_j^2(t)=\int_0^t \sigma_{\Lambda_j}^{-1}(h_{(j,t)})(u)\, d\Lambda_{j,0}(u).
\end{eqnarray*}
\end{theorem}

We now turn to the issue of estimating the asymptotic variances of the estimators $\widehat\beta_n$, $\widehat\gamma_n$, and $\widehat\Lambda_{j,n}(t)$ ($t\in(0,\tau)$, $j\in\mathcal K$). It turns out that the asymptotic variances $\Sigma_\beta$, $\Sigma_\gamma$, and $v_j^2(t)$ are not expressible in explicit forms, since the inverse $\sigma^{-1}$ has no closed form. However, this is not a problem if we can provide consistent estimators for them. Such estimators are defined below.

For $i=1,\ldots,n$, let $X_{ir}$ denote the $r$-th ($r=1,\ldots,p$) component of $X_i$, $S_{\gamma,i}(\theta)$ be defined as in (\ref{score}) with $\mathcal O$ and $W$ replaced by $\mathcal O_i$ and $W_i$ respectively, and $S_{\gamma,i,s}(\theta)$ be the $s$-th ($s=1,\ldots,q$) component of $S_{\gamma,i}(\theta)$. Using these notations, we define the following block matrix
\begin{eqnarray}\label{bm}
\mathbb A_n=\left(
\begin{array}{ccc}
A^{\beta\beta}& A^{\beta\gamma}& A^{\beta\Lambda}\\
A^{\gamma\beta}& A^{\gamma\gamma}& A^{\gamma\Lambda}\\
A^{\Lambda\beta}& A^{\Lambda\gamma}& A^{\Lambda\Lambda}\\
\end{array}
\right)
\end{eqnarray}
where the sub-matrices $A^{\beta\beta}, A^{\gamma\gamma}, A^{\beta\gamma}$, and $A^{\gamma\beta}$ are defined as follows by their $(r,s)$-th component:
\begin{eqnarray*}
&&A^{\beta\beta}_{rs}=\frac{1}{n}\sum_{i=1}^n\{\psi(\mathcal O_i,\widehat\theta_n)\}^2X_{ir}X_{is}, \quad r,s=1,\ldots,p,\\
&&A^{\gamma\gamma}_{rs}=\frac{1}{n}\sum_{i=1}^nS_{\gamma,i,r}(\widehat\theta_n)S_{\gamma,i,s}(\widehat\theta_n), \quad r,s=1,\ldots,q,\\
&&A^{\beta\gamma}_{rs}=\frac{1}{n}\sum_{i=1}^n\psi(\mathcal O_i,\widehat\theta_n) X_{ir} S_{\gamma,i,s}(\widehat\theta_n), \quad r=1,\ldots,p, \quad s=1,\ldots,q,
\\
&&A^{\gamma\beta}_{rs}=A^{\beta\gamma}_{sr}, \quad r=1,\ldots,q, \quad s=1,\ldots,p.
\end{eqnarray*}
Define the block matrices $A^{\beta\Lambda}=(A^{\beta\Lambda_1},\ldots,A^{\beta\Lambda_K})$ and $A^{\gamma\Lambda}=(A^{\gamma\Lambda_1},\ldots,A^{\gamma\Lambda_K})$, where for every $k\in\mathcal K$, the sub-matrices $A^{\beta\Lambda_k}$ and $A^{\gamma\Lambda_k}$ are defined by
\begin{eqnarray*}
&&A^{\beta\Lambda_k}_{rs}=\frac{2}{n}X_{sr}\Delta_s\psi(\mathcal O_s,\widehat\theta_n)Q(\mathcal O_s,k,\widehat\theta_n), \quad r=1,\ldots,p, \quad s=1,\ldots,n,\\
&&A^{\gamma\Lambda_k}_{rs}=\frac{2}{n}S_{\gamma,s,r}(\widehat\theta_n)\Delta_sQ(\mathcal O_s,k,\widehat\theta_n), \quad r=1,\ldots,q, \quad s=1,\ldots,n.
\end{eqnarray*}
Define also the block matrices
\begin{eqnarray*}
A^{\Lambda\beta}=\left(
\begin{array}{c}
A^{\Lambda_1\beta}\\
\vdots\\
A^{\Lambda_K\beta}
\end{array}
\right)\quad\quad
A^{\Lambda\gamma}=\left(
\begin{array}{c}
A^{\Lambda_1\gamma}\\
\vdots\\
A^{\Lambda_K\gamma}
\end{array}
\right)\quad\quad
A^{\Lambda\Lambda}=\left(
\begin{array}{ccc}
A^{\Lambda_1\Lambda_1}& \ldots& A^{\Lambda_1\Lambda_K}\\
\vdots& &\vdots\\
A^{\Lambda_K\Lambda_1}& \ldots& A^{\Lambda_K\Lambda_K}
\end{array}
\right)
\end{eqnarray*}
where for every $j,k\in\mathcal K$,
\begin{eqnarray*}
&&A^{\Lambda_k\beta}_{rs}=-\frac{1}{n}\sum_{i=1}^n2X_{is}\psi(\mathcal O_i,\widehat\theta_n)Q(\mathcal O_i,k,\widehat\theta_n) e^{\widehat\beta_n'X_i} Y_i(T_r), \quad r=1,\ldots,n, \quad s=1,\ldots,p,\\
&&A^{\Lambda_k\gamma}_{rs}=-\frac{1}{n}\sum_{i=1}^n2S_{\gamma,i,s}(\widehat\theta_n)Q(\mathcal O_i,k,\widehat\theta_n) e^{\widehat\beta_n'X_i} Y_i(T_r), \quad r=1,\ldots,n, \quad s=1,\ldots,q,\\
&&A^{\Lambda_k\Lambda_j}_{rs}=1\{j=k\}1\{r=s\} \frac{1}{n}\sum_{i=1}^n Q(\mathcal O_i,k,\widehat\theta_n)\phi(T_r, \mathcal O_i,k,\widehat\theta_n)\\
&&\hspace{1.25cm}+1\{j>k\}\left(1\{r=s\}\frac{1}{n}\sum_{i=1}^n 2 \phi(T_s, \mathcal O_i,k,\widehat\theta_n)Q(\mathcal O_i,j, \widehat\theta_n)\right.\\
&&\hspace{1.5cm} +\frac{2}{n} \sum_{i=1}^n\phi(T_r, \mathcal O_i, k, \widehat\theta_n)Q(\mathcal O_i,j, \widehat\theta_n)e^{\widehat\beta_n'X_i} \widehat{\Delta \Lambda_{j,n}}(T_s) \{1\{T_s\leq T_i\}- 1\{T_s\leq T_r\}\}\\
&&\hspace{1.5cm} \left.-\frac{2}{n}\phi(T_r, \mathcal O_s,k,\widehat\theta_n) Q(\mathcal O_s,j, \widehat\theta_n)\Delta_s\right), \quad r,s=1,\ldots,n,
\end{eqnarray*}
and $\widehat{\Delta \Lambda_{j,n}}(T_s)$ is the jump size of $\widehat\Lambda_{j,n}$ at $T_s$ that is, $\widehat{\Delta \Lambda_{j,n}}(T_s)=\widehat\Lambda_{j,n}(T_s)-\widehat\Lambda_{j,n}(T_s-)$ ($j\in\mathcal K, s=1,\ldots,n$). Note that for notational simplicity, the lower (sample size) indice $n$ has been omitted in the notations for the sub-matrices of $\mathbb A_n$.

Now, define 
\begin{eqnarray*}
&&\widehat\Sigma_{\beta,n}=\left\{A^{\beta\beta}-A^{\beta\gamma}(A^{\gamma\gamma})^{-1}A^{\gamma\beta}-(A^{\beta\Lambda}- A^{\beta\gamma}(A^{\gamma\gamma})^{-1}A^{\gamma\Lambda})\right.\\
&&\hspace{1.5cm}\left.\times(A^{\Lambda\Lambda}-A^{\Lambda\gamma}(A^{\gamma\gamma})^{-1}A^{\gamma\Lambda})^{-1}(A^{\Lambda\beta}-A^{\Lambda\gamma}(A^{\gamma\gamma})^{-1}A^{\gamma\beta})\right\}^{-1},\\
&&\widehat\Sigma_{\gamma,n}=\left\{A^{\gamma\gamma}-A^{\gamma\beta}(A^{\beta\beta})^{-1}A^{\beta\gamma}-(A^{\gamma\Lambda}- A^{\gamma\beta}(A^{\beta\beta})^{-1}A^{\beta\Lambda})\right.\\
&&\hspace{1.5cm}\left.\times(A^{\Lambda\Lambda}-A^{\Lambda\beta}(A^{\beta\beta})^{-1}A^{\beta\Lambda})^{-1}(A^{\Lambda\gamma}-A^{\Lambda\beta}(A^{\beta\beta})^{-1}A^{\beta\gamma})\right\}^{-1},
\end{eqnarray*}
and
\begin{eqnarray*}
&&\widehat\Sigma_{\Lambda,n} =\left\{A^{\Lambda\Lambda}-A^{\Lambda\beta}(A^{\beta\beta})^{-1}A^{\beta\Lambda}-(A^{\Lambda\gamma}- A^{\Lambda\beta}(A^{\beta\beta})^{-1}A^{\beta\gamma})\right.\\
&&\hspace{1.5cm}\left.\times(A^{\gamma\gamma}-A^{\gamma\beta}(A^{\beta\beta})^{-1}A^{\beta\gamma})^{-1}(A^{\gamma\Lambda}-A^{\gamma\beta}(A^{\beta\beta})^{-1}A^{\beta\Lambda})\right\}^{-1}.
\end{eqnarray*}

Then the following holds:

\begin{theorem}\label{varasymp}
Under conditions (a)-(f), $\widehat\Sigma_{\beta,n}$ and $\widehat\Sigma_{\gamma,n}$ converge in probability to $\Sigma_\beta$ and $\Sigma_\gamma$ respectively as $n$ tends to $\infty$. Moreover, for $t\in(0,\tau)$ and $j\in\mathcal K$, let
\begin{eqnarray*}
\widehat v^2_{j,n}(t)=\widehat \Xi_{(j,t)}^{n'}\widehat\Sigma_{\Lambda,n} U_{(j,t)}^n,
\end{eqnarray*}
where $$\widehat \Xi_{(j,t)}^n=\left(0_{(j-1)n}',\widehat{\Delta \Lambda_{j,n}}(T_1) 1\{T_1\leq t\},\ldots,\widehat{\Delta \Lambda_{j,n}}(T_n) 1\{T_n\leq t\},0_{(K-j)n}'\right)'$$ and $$U_{(j,t)}^n=(0_{(j-1)n}',1\{T_1\leq t\},\ldots,1\{T_n\leq t\},0_{(K-j)n}')'.$$ Then $\widehat v^2_{j,n}(t)$ converges  in probability to $v^2_j(t)$ as $n$ tends to $\infty$.
\end{theorem}

\section{Discussion}\label{disc}

In this paper, we have constructed consistent and asymptotically normal estimators for the stratified proportional intensity regression model when the sample stratum information is only partially available. The proposed estimator for the regression parameter of interest in this model has been shown to be semiparametrically efficient. Although computationally more challenging, these estimators improve the ones previously investigated in the literature, such as the regression calibration estimators \citep{duplec}.

We have obtained explicit (and computationally fairly simple) formulas for consistent estimators of the asymptotic variances. These formulas may however require the inversion of potentially large matrices. For a large sample, this inversion may be unstable. An alternative solution relies on numerical differentiation of the profile log-likelihood (see \cite{mrvdv} and \cite{chen99} for example). Note that in this latter method however, no estimator is available for the asymptotic variance of the cumulative baseline intensity estimator. Some further work is needed to evaluate the numerical performance of the proposed estimators. This is the subject for future research, and requires some extensive simulation work which falls beyond the scope of this paper.

In this paper, a multinomial logistic model \citep{jobson} is used for modeling the conditional stratum probabilities given covariates. This choice was mainly motivated by the fact that this model is commonly used in medical research for modeling the relationship between a categorical response and covariates. The theoretical results developed here can be extended to the case of other link functions. In addition, the covariate $X$ in model (\ref{model}) is assumed to be time independent, for convenience. This assumption can be relaxed to accomodate time varying covariates, provided that appropriate regularity conditions are made.

\section*{Appendix A. Proofs of Theorems}

\textbf{A.1 Proof of Theorem \ref{exist}}
\\\\
For every $k\in\mathcal K$, define $\mathcal I_{k}^n=\{i\in\{1,\ldots,n\}|\Delta_iR_i1\{S_i=k\}=1 \mbox{ or } \Delta_i(1-R_i)=1\}$, and let $i_{k}^n$ denote the cardinality of $\mathcal I_{k}^n$. Let $i_{\bullet}^n=\sum_{k=1}^Ki_{k}^n$. Consider the set of times $\{T_i,i\in\mathcal I_k^n\}$. Let $t_{(k,1)}<\ldots<t_{(k,i^n_{k})}$ denote the ordered failure times in this set. For any given sample size $n$, the NPML estimation method consists in maximizing $L_n$ in (\ref{odl}) over the approximating parameter space 
\begin{eqnarray*}
\Theta_n=\left\{(\beta,\gamma,\Lambda_k\{t_{(k,j)}\}):\beta\in\mathcal B; \gamma\in\mathcal G; \Lambda_k\{t_{(k,j)}\}\in[0,\infty), j=1,\ldots,i_{k}^n, k\in\mathcal K\right\}.
\end{eqnarray*}
Suppose first that $\Lambda_k\{t_{(k,j)}\}\leq M<\infty$, for $j=1,\ldots,i_{k}^n$ and $k\in\mathcal K$. Since $L_n$ is a continuous function of $\beta,\gamma$, and the $\Lambda_k\{t_{(k,j)}\}$'s on the compact  set $\mathcal B\times\mathcal G\times [0,M]^{i_{\bullet}^n}$, $L_n$ achieves its maximum on this set.
\\
To show that a maximum of $L_n$ exists on $\mathcal B\times\mathcal G\times [0,\infty)^{i_{\bullet}^n}$, we show that there exists a finite $M$ such that for all $\theta^M=(\beta^M,\gamma^M,(\Lambda_k^M\{t_{(k,j)}\})_{j,k})\in(\mathcal B\times\mathcal G\times [0,\infty)^{i_{\bullet}^n})\backslash (\mathcal B\times\mathcal G\times [0,M]^{i_{\bullet}^n})$, there exists a $\theta=(\beta,\gamma,(\Lambda_k\{t_{(k,j)}\})_{j,k})\in\mathcal B\times\mathcal G\times [0,M]^{i_{\bullet}^n}$ such that $L_n(\theta)> L_n(\theta^M)$. A proof by contradiction is adopted for that purpose.
\\
Assume that for all $M<\infty$, there exists $\theta^M\in(\mathcal B\times\mathcal G\times [0,\infty)^{i_{\bullet}^n}) \backslash (\mathcal B\times\mathcal G\times [0,M]^{i_{\bullet}^n})$ such that for all $\theta\in \mathcal B\times\mathcal G\times [0,M]^{i_{\bullet}^n}$, $L_n(\theta)\leq L_n(\theta^M)$. It can be seen that $L_n$ is bounded above by
\begin{eqnarray*}
K^n\prod_{i=1}^n \left[\prod_{k=1}^K\left\{c_3\Lambda_k\{T_i\}\right\}^{\Delta_iR_i1\{S_i=k\}}\exp\left(-c_2R_i1\{S_i=k\}\sum_{j=1}^{i_{k}^n}\Lambda_k\{t_{(k,j)}\}1\{t_{(k,j)}\leq T_i\}\right)\right].
\end{eqnarray*}
If $\theta^M\in (\mathcal B\times\mathcal G\times [0,\infty)^{i_{\bullet}^n}) \backslash (\mathcal B\times\mathcal G\times [0,M]^{i_{\bullet}^n})$, then there exists $l\in\mathcal K$ and $p\in\{1,\ldots,i_l^n\}$ such that $\Lambda_l^M\{t_{(l,p)}\}>M$.
By assumption (d), there exists at least one individual with indice $i_M$ $(i_M\in\{1,\ldots,n\})$ such that $1\{S_{i_M}=l\}=1$, $Y_{i_M}(\tau)=1$ (and therefore $t_{(l,p)}\leq T_{i_M}=\tau$), and $R_{i_M}=1$. Hence 
\begin{eqnarray*}
R_{i_M}1\{S_{i_M}=l\}\sum_{j=1}^{i_l^n}\Lambda^M_l\{t_{(l,j)}\}1\{t_{(l,j)}\leq T_{i_M}\}\rightarrow \infty \mbox{ as } M\rightarrow \infty.
\end{eqnarray*}
It follows that the upper bound of $L_n(\theta^M)$ (and therefore $L_n(\theta^M)$ itself) can be made as close to 0 as desired by increasing $M$. This is the desired contradiction.
\\
\\
$\Box$
\\
\\
\textbf{A.2 Proof of Theorem \ref{csy}}
\\\\
We adapt the techniques developed by \cite{murphy94}, in order to prove consistency of our proposed estimator $\widehat\theta_n$. The proof essentially consists of three steps: (i) for every $k\in\mathcal K$, we show that the sequence $\widehat\Lambda_{k,n}(\tau)$ is almost surely bounded as $n$ goes to infinity, (ii) we show that every subsequence of $n$ contains a further subsequence along which the NPMLE $\widehat\theta_n$ converges, (iii) we show that the limit of every convergent subsequence of $\widehat\theta_n$ is $\theta_0$.
\\\\
\textit{Proof of (i).} Note first that for all $s\in[0,\tau]$ and $k\in\mathcal K$, $\frac{1}{n}\sum_{i=1}^nQ(\mathcal O_i,k,\widehat\theta_n) e^{\widehat\beta_n'X_i}Y_i(s)\geq c_2\frac{1}{n}\sum_{i=1}^nR_i1\{S_i=k\}Y_i(\tau)$. Moreover, $Q(\mathcal O_i,k,\widehat\theta_n)$ is bounded by 1. It follows that for all $k\in\mathcal K$,
\begin{eqnarray*}
0\leq \widehat\Lambda_{k,n}(\tau)\leq \frac{1}{c_2}\int_0^\tau\frac{d\bar N_n(s)}{\frac{1}{n}\sum_{i=1}^nR_i1\{S_i=k\}Y_i(\tau)}= \frac{\frac{1}{n}\sum_{i=1}^n \Delta_i}{c_2\frac{1}{n}\sum_{i=1}^nR_i1\{S_i=k\}Y_i(\tau)},
\end{eqnarray*}
where $\bar N_n(s)=n^{-1}\sum_{i=1}^n N_i(s)$. Next, $\frac{1}{n}\sum_{i=1}^nR_i1\{S_i=k\}Y_i(\tau)$ converges almost surely to $P_{\theta_0}[R1\{S=k\}Y(\tau)]>c_4>0$ therefore, for each $k\in\mathcal K$, as $n$ goes to infinity, $\widehat\Lambda_{k,n}(\tau)$ is bounded above almost surely by $\frac{1}{c_2c_4}$.
\\\\
\textit{Proof of (ii).} If (i) holds, by Helly's theorem (see \cite{loeve}, p179), every subsequence of $n$ has a further subsequence along which $\widehat\Lambda_{1,n}$ converges weakly to some nondecreasing right-continuous function $\Lambda_1^\ast$, with probability 1. By successive extractions of sub-subsequences, we can further find a subsequence (say $n_j$) such that $\widehat\Lambda_{k,n_j}$ converges weakly to some nondecreasing right-continuous function $\Lambda_k^\ast$, for every $k\in\mathcal K$, with probability 1. By the compactness of $\mathcal B\times\mathcal G$, we can further find a subsequence of $n_j$ (we shall still denote it by $n_j$ for simplicity of notations) such that $\widehat\Lambda_{k,n_j}$ converges weakly to $\Lambda_k^\ast$ (for every $k\in\mathcal K$) and $(\widehat\beta_{n_j}, \widehat\gamma_{n_j})$ converges to some $(\beta^\ast,\gamma^\ast)$, with probability 1. We now show that the $\Lambda_k^\ast$'s must be continuous on $[0,\tau]$.

Let $\psi$ be any nonnegative, bounded, continuous function. Then, for any given $k\in\mathcal K$,
\begin{eqnarray*}
&&\int^\tau_0\psi(s)\,d\Lambda^\ast_k(s)=\int^\tau_0\psi(s)\,d\{\Lambda^\ast_k(s)-\widehat\Lambda_{k,n_j}(s)\}\\
&&\hspace{.5cm}+\int^\tau_0\psi(s)\left[\frac{1}{n_j}\sum_{l=1}^{n_j}Q(\mathcal O_l,k,\widehat\theta_{n_j})e^{\widehat\beta_{n_j}'X_l}Y_l(s)\right]^{-1}\frac{1}{n_j}\sum_{i=1}^{n_j}Q(\mathcal O_i,k,\widehat\theta_{n_j})\,dN_i(s) \\
&&\hspace{.5cm}\leq\int^\tau_0\psi(s)\,d\{\Lambda^\ast_k(s)-\widehat\Lambda_{k,n_j}(s)\}+\int^\tau_0\psi(s)\left[\frac{c_2}{n_j}\sum_{l=1}^{n_j}R_l1\{S_l=k\}Y_l(s)\right]^{-1}d\bar N_{n_j}(s).
\end{eqnarray*}
 By the Helly-Bray Lemma (see \cite{loeve}, p180), $\int^\tau_0\psi(s)\,d\{\Lambda^\ast_k(s)-\widehat\Lambda_{k,n_j}(s)\}\rightarrow 0$ as $j\rightarrow \infty$. Moreover, $\bar N_{n_j}(\cdot)$ and $\frac{1}{n_j}\sum_{l=1}^{n_j}R_l1\{S_l=k\}Y_l(\cdot)$ converge almost surely in supremum norm to
\begin{eqnarray*}
\sum_{k=1}^K\int_0^\cdot P_{\theta_0}\left[1\{S=k\}e^{\beta_0'X}Y(s)\right]\,d\Lambda_{k,0}(s)\mbox{ and } P_{\theta_0}\left[R1\{S=k\}Y(\cdot)\right]
\end{eqnarray*}
respectively, where the latter term is bounded away from 0 on $s\in[0,\tau]$ by assumption (d). Thus, by applying the extended version of the Helly-Bray Lemma (stated by \cite{kor} for example) to the second term on the right-hand side of the previous inequality, we get that
\begin{eqnarray}\label{ineg}
&&\int^\tau_0\psi(s)\,d\Lambda^\ast_k(s)\\
&&\hspace{1cm}\leq c_2\int^\tau_0\psi(s)\left\{P_{\theta_0}\left[R1\{S=k\}Y(s)\right]\right\}^{-1}\sum_{k=1}^KP_{\theta_0}[1\{S=k\}e^{\beta_0'X}Y(s)]\lambda_{k,0}(s)\,ds.\nonumber\\
&&\hspace{1cm}\leq \frac{c_2c_3}{c_4} \sum_{k=1}^K \int^\tau_0\psi(s) \lambda_{k,0}(s)\,ds.\nonumber
\end{eqnarray}
Suppose that $\Lambda^\ast_k$ has discontinuities, and let $\psi$ be close to 0 except at the jump points of $\Lambda^\ast_k$, where it is allowed to have high and thin peaks. While the right-hand side of inequality (\ref{ineg}) should be close to 0 ($\lambda_{k,0}$ is continuous by assumption (b)), its left-hand side can be made arbitrarily large, yielding a contradiction. Thus $\Lambda^\ast_k$ must be continuous ($k\in\mathcal K$). A second conclusion, arising from Dini's theorem, is that $\widehat\Lambda_{k,n_j}$ uniformly converges to $\Lambda^\ast_k$ ($k\in\mathcal K$), with probability 1. To summarize: for any given subsequence of $n$, we have found a further subsequence $n_j$ and an element $(\beta^\ast,\gamma^\ast, \Lambda^\ast_k, k\in\mathcal K)$ such that $\|\widehat\beta_{n_j}-\beta^\ast\|$, $\|\widehat\gamma_{n_j}-\gamma^\ast\|$, and $\sup_{t\in[0,\tau]}|\widehat\Lambda_{k,n_j}(t)-\Lambda^\ast_k(t)|$ (for every $k\in\mathcal K$) converge to 0 almost surely.
\\\\
\textit{Proof of (iii).}
To prove (iii), we first define random step functions
\begin{eqnarray*}
\overline\Lambda_{k,n}(t)=\int_0^t\sum_{i=1}^n \frac{Q(\mathcal O_i,k,\theta_0)}{\sum_{j=1}^n Q(\mathcal O_j,k,\theta_0)\exp(\beta_0'X_j)Y_j(s)}\,dN_i(s),\; 0\leq t\leq \tau, k\in\mathcal K,
\end{eqnarray*}
and we show that for every $k\in\mathcal K$, $\overline\Lambda_{k,n}$ almost surely uniformly converges to $\Lambda_{k,0}$ on $[0,\tau]$. First, note that
\begin{eqnarray}\label{ineq}
&&\sup_{t\in [0,\tau]}\left|\overline\Lambda_{k,n}(t)-P_{\theta_0}\left[\frac{\Delta 1\{T\leq t\}Q(\mathcal O,k,\theta_0)}{P_{\theta_0}\left[1\{S=k\}e^{\beta_0'X}Y(s)\right]\left.\right|_{s=T}}\right]\right|\nonumber\\
&&\hspace{1cm}\leq\sup_{t\in [0,\tau]}\left|\frac{1}{n}\sum_{i=1}^n \Delta_i1\{T_i\leq t\}Q(\mathcal O_i,k,\theta_0)\right.\nonumber\\
&&\left.\hspace{2cm}\times\left\{\frac{1}{\mathbb P_n\left[Q(\mathcal O,k,\theta_0)e^{\beta_0'X}Y(s)\right]}-\frac{1}{P_{\theta_0}\left[1\{S=k\}e^{\beta_0'X}Y(s)\right]}\right\}_{\left.\right|_{s=T_i}}\right|\nonumber\\
&&\hspace{3cm}+\sup_{t\in [0,\tau]}\left|(\mathbb P_n-P_{\theta_0})\left[\frac{\Delta 1\{T\leq t\}Q(\mathcal O,k,\theta_0)}{P_{\theta_0}\left[1\{S=k\}e^{\beta_0'X}Y(s)\right]\left.\right|_{s=T}}\right]\right|\nonumber\\
&&\hspace{1cm}\leq\sup_{s\in [0,\tau]}\left|\frac{1}{\mathbb P_n\left[Q(\mathcal O,k,\theta_0)e^{\beta_0'X}Y(s)\right]}-\frac{1}{P_{\theta_0}\left[1\{S=k\}e^{\beta_0'X}Y(s)\right]}\right|\nonumber\\
&&\hspace{3cm}+\sup_{t\in [0,\tau]}\left|(\mathbb P_n-P_{\theta_0})\left[\frac{\Delta 1\{T\leq t\}Q(\mathcal O,k,\theta_0)}{P_{\theta_0}\left[1\{S=k\}e^{\beta_0'X}Y(s)\right]\left.\right|_{s=T}}\right]\right|
\end{eqnarray}
The class $\{Y(s):s\in [0,\tau]\}$ is Donsker and $Q(\mathcal O,k,\theta_0)e^{\beta_0'X}$ is a bounded measurable function, hence $\{Q(\mathcal O,k,\theta_0)e^{\beta_0'X}Y(s):s\in [0,\tau]\}$ is Donsker (Corollary 9.31, \cite{k07}), and therefore Glivenko-Cantelli. Moreover, $P_{\theta_0}[Q(\mathcal O,k,\theta_0)e^{\beta_0'X}Y(s)]=P_{\theta_0}[P_{\theta_0}[1\{S=k\}|\mathcal O] e^{\beta_0'X}Y(s)]=P_{\theta_0}[1\{S=k\}e^{\beta_0'X}Y(s)]$. Thus
\begin{eqnarray*}
\sup_{s\in [0,\tau]}\left|\mathbb P_n\left[Q(\mathcal O,k,\theta_0)e^{\beta_0'X}Y(s)\right]-P_{\theta_0}\left[1\{S=k\}e^{\beta_0'X}Y(s)\right]\right|
\end{eqnarray*}
converges to 0 a.e. Next, $P_{\theta_0}[1\{S=k\}e^{\beta_0'X}Y(s)]$ is larger than $c_2.P_{\theta_0}[1\{S=k\}Y(\tau)]$ and thus, by assumption (d), $P_{\theta_0}[1\{S=k\}e^{\beta_0'X}Y(s)]>0$. It follows that the first term on the right-hand side of inequality (\ref{ineq}) converges to 0 a.e.. Similar aguments show that the class $\{\Delta 1\{T\leq t\}Q(\mathcal O,k,\theta_0)\slash P_{\theta_0}[1\{S=k\}e^{\beta_0'X}Y(s)]\left.\right|_{s=T}:t\in [0,\tau]\}$ is also a Glivenko-Cantelli class, and therefore $\overline\Lambda_{k,n}$ almost surely uniformly converges to 
\begin{eqnarray*}
P_{\theta_0}\left[\frac{\Delta 1\{T\leq t\}Q(\mathcal O,k,\theta_0)}{P_{\theta_0}\left[1\{S=k\}e^{\beta_0'X}Y(s)\right]\left.\right|_{s=T}}\right].
\end{eqnarray*}
Now, note that $\Lambda_{k,0}(t)=\int_0^t\frac{P_{\theta_0}[1\{S=k\}\,dN(s)]}{P_{\theta_0}[1\{S=k\}e^{\beta_0'X}Y(s)]}$, which can be reexpressed as
\begin{eqnarray*}
\Lambda_{k,0}(t)=\frac{P_{\theta_0}\left[1\{S=k\}\Delta1\{T\leq t\}\right]}{P_{\theta_0}\left[1\{S=k\}e^{\beta_0'X}Y(s)\right]\left.\right|_{s=T}}=P_{\theta_0}\left[\frac{\Delta 1\{T\leq t\}Q(\mathcal O,k,\theta_0)}{P_{\theta_0}\left[1\{S=k\}e^{\beta_0'X}Y(s)\right]\left.\right|_{s=T}}\right].
\end{eqnarray*}
Thus $\overline\Lambda_{k,n}$ almost surely uniformly converges to $\Lambda_{k,0}$ on $[0,\tau]$.

Next, using somewhat standard arguments (see \cite{parner} for example), we can show that $0\leq n_j^{-1}\{\log L_{n_j}(\widehat\theta_{n_j})-\log L_{n_j}(\overline\theta_{n_j})\}$ converges to the negative Kullback-Leibler information $P_{\theta_0}[\log(L_1(\theta^\ast)\slash L_1(\theta_0))]$. Thus, the Kullback-Leibler information must be zero, and it follows that with probability 1, $L_1(\theta^\ast)=L_1(\theta_0)$. The proof of consistency is completed if we show that this equality implies $\theta^\ast=\theta_0$. For that purpose, consider $L_1(\theta^\ast)=L_1(\theta_0)$ under $\Delta=1$, $R=1$, and $1\{S=k\}=1$ (for each $k\in \mathcal K$ in turn). Note that this is possible by assumption (e). This yields the following equation for almost all $t\in[0,\tau], \|x\|<c_1, \|w\|<c_1$:
\begin{eqnarray*}
\log\frac{\lambda_k^\ast(t)}{\lambda_{k,0}(t)}+\left(\beta^\ast-\beta_0\right)'x-\Lambda_k^\ast(t)e^{\beta^{\ast'}x}+\Lambda_{k,0}(t)e^{\beta_0'x}+\log\frac{\pi^\ast_k(w)}{\pi_{k,0}(w)}=0.
\end{eqnarray*}
This equation is analogous to equation (A.2) in \cite{chen99}. The rest of the proof of identifiability thus proceeds along the same lines as the proof of Lemma A.1.1 in \cite{chen99}, and is omitted.
\\
 Hence, for any given subsequence of $n$, we have found a further subsequence $n_j$ such that $\|\widehat\beta_{n_j}-\beta_0\|$, $\|\widehat\gamma_{n_j}-\gamma_0\|$, and $\sup_{t\in[0,\tau]}|\widehat\Lambda_{k,n_j}(t)-\Lambda_{k,0}(t)|$ (for every $k\in\mathcal K$) converge to 0 almost surely, which implies that the sequence of NPMLE $\widehat\theta_n$ converges almost surely to $\theta_0$.
\\
\\
$\Box$
\\
\\
\textbf{A.3 Proof of Theorem \ref{asympn}}
\\
\\
The proof of Theorem \ref{asympn} uses similar arguments as the proof of Theorem 3 of \cite{fang05}, so we only highlight the parts that are different. We need a few lemmas before presenting the proof.
\begin{lemma}\label{espsc0}
Let $h\in H$. Then the following holds: $P_{\theta_0}\left[S_1(\theta_0)(h)\right]=P_{\theta_0}[h_\beta' S_\beta(\theta_0)+ h_\gamma'S_\gamma(\theta_0)+\sum_{k=1}^K S_{\Lambda_k}(\theta_0)(h_{\Lambda_k})]=0$.
\end{lemma}

\textbf{Proof.} From the properties of the conditional expectation, we first note that
\begin{eqnarray*}
P_{\theta_0}\left[S_\beta(\theta_0)\right]&=&P_{\theta_0}\left[\Delta X-\sum_{k=1}^K Q(\mathcal O,k,\theta_0)X\exp(\beta_0'X)\Lambda_{k,0}(T)\right]\\
&=&P_{\theta_0}\left[\Delta X-\sum_{k=1}^K 1\{S=k\}X\exp(\beta_0'X)\Lambda_{k,0}(T)\right]\\
&=&P_{\theta_0}\left[XM(\tau)\right],
\end{eqnarray*}
where $M(t)=N(t)-\int_0^t\sum_{k=1}^K1\{S=k\}e^{\beta_0'X}Y(u)\,d\Lambda_{k,0}(u)$ is the counting process martingale with respect to the filtration $\mathcal F_t=\sigma\{N(u), 1\{C\leq u\}, X, S, W:0\leq u\leq t\}$. $X$ is bounded and $\mathcal F_t$-measurable, hence it follows that $P_{\theta_0}[S_\beta(\theta_0)]=0$. Using similar arguments, we can verify that $P_{\theta_0}[S_{\Lambda_k}(\theta_0)(h_{\Lambda_k})]=0$, $k\in\mathcal K$. Finally, for $k=1,\ldots,K-1$,
\begin{eqnarray*}
P_{\theta_0}\left[S_{\gamma_k}(\theta_0)\right]&=&P_{\theta_0}\left[W\left[Q(\mathcal O,k,\theta_0)-\pi_{k,\gamma_0}(W)\right]\right]\\
&=&P_{\theta_0}\left[WP_{\theta_0}\left[1\{S=k\}-\pi_{k,\gamma_0}(W)|W\right]\right]\\
&=&0.
\end{eqnarray*}
Combining these results yields that $P_{\theta_0}\left[S_1(\theta_0)(h)\right]=0$.
\\
\\
$\Box$
\\
\\
We now come to the continuous invertibility of the continuous linear operator $\sigma$ defined in Section \ref{asymprp}.
\begin{lemma}\label{infoope}
The operator $\sigma$ is continuously invertible.
\end{lemma}

\textbf{Proof.} Since $H$ is a Banach space, to prove that $\sigma$ is continuously invertible, it is sufficient to prove that $\sigma$ is one-to-one and that it can be written as the sum of a bounded linear operator with a bounded inverse and a compact operator (Lemma 25.93 of \cite{vdv98}).

Define the linear operator $A(h)=(h_\beta, h_\gamma, P_{\theta_0}\left[ 1\{S=k\} \phi(\cdot,\mathcal O,k,\theta_0)\right]h_{\Lambda_k}(\cdot);k\in\mathcal K)$, this is a bounded operator due to the boundedness of $X$. Moreover, for all $u\in[0,\tau]$ and $k\in\mathcal K$, $P_{\theta_0}\left[ 1\{S=k\} \phi(u,\mathcal O,k,\theta_0)\right]\geq c_2c_4>0$ by assumptions (c) and (d). This implies that $A$ is invertible with bounded inverse $A^{-1}(h)=(h_\beta, h_\gamma, P_{\theta_0}\left[ 1\{S=k\} \phi(\cdot,\mathcal O,k,\theta_0)\right]^{-1}h_{\Lambda_k}(\cdot);k\in\mathcal K)$. The operator $\sigma-A$ can be shown to be compact by using the same techniques as in \cite{lu08} for example.

To prove that $\sigma$ is one-to-one, let $h\in H$ such that $\sigma(h)=0$. If $\sigma(h)=0$, $P_{\theta_0}\left[S_1(\theta_0)(h)^2\right]=0$, and therefore $S_1(\theta_0)(h)=0$ almost surely. Let $j\in\mathcal K$. By assumption (e), for almost every $t\in[0,\tau], \|x\|\leq c_1$, and $\|w\|\leq c_1$, there is a non-negligible set $\Omega_{t,x,w}\subseteq\Omega$ such that $\Delta(\omega)=1$, $R(\omega)=1$, and $1\{S(\omega)=j\}=1$ when $\omega\in\Omega_{t,x,w}$. If $S_1(\theta_0)(h)=0$ almost surely, then in particular, for almost every $t\in[0,\tau], \|x\|\leq c_1$, and $\|w\|\leq c_1$, $S_1(\theta_0)(h)=0$ when $\omega\in\Omega_{t,x,w}$, which yields the following equation:
\begin{eqnarray}\label{inj1}
h_{\Lambda_j}(t)+h_\beta'x+w'h_{\gamma_j}-\sum_{k=1}^{K-1}w'h_{\gamma_k}\pi_{k,\gamma_0}(w)-e^{\beta_0'x}\left[\int_0^t h_{\Lambda_j}(s)\,d\Lambda_{j,0}(s)+h_\beta'x \Lambda_{j,0}(t)\right]=0,
\end{eqnarray}
with $h_{\gamma_j}=0$ when $j=K$. Then, by choosing $t$ arbitrarily close to 0, and since $\Lambda_{j,0}$ is continuous, $\Lambda_{j,0}(0)=0$, and $h_{\Lambda_j}$ is continuous from the right at 0, we get that
\begin{eqnarray}\label{inj2}
h_{\Lambda_j}(0)+h_\beta'x+w'h_{\gamma_j}-\sum_{k=1}^{K-1}w'h_{\gamma_k}\pi_{k,\gamma_0}(w)=0.
\end{eqnarray}
Taking the difference (\ref{inj1})-(\ref{inj2}) yields that
\begin{eqnarray}\label{inj3}
h_{\Lambda_j}(t)-h_{\Lambda_j}(0)=e^{\beta_0'x}\left[\int_0^t h_{\Lambda_j}(s)\,d\Lambda_{j,0}(s)+h_\beta'x \Lambda_{j,0}(t)\right]
\end{eqnarray}
for almost every $t\in[0,\tau]$ and $\|x\|\leq c_1$. Since $\Lambda_{j,0}$ is increasing (by assumption (b)), for every $t>0$, $\Lambda_{j,0}(t)>\Lambda_{j,0}(0)=0$ and therefore (\ref{inj3}) can be rewritten as
\begin{eqnarray}\label{inj4}
\frac{h_{\Lambda_j}(t)-h_{\Lambda_j}(0)}{\Lambda_{j,0}(t)}=e^{\beta_0'x}\left[r(t)+h_\beta'x \right],
\end{eqnarray}
where $r(t)=\int_0^t h_{\Lambda_j}(s)\,d\Lambda_{j,0}(s)\slash\Lambda_{j,0}(t)$. Consider first the case where $\beta_0=0$. Since the left-hand side of (\ref{inj4}) does not depend on $x$, $h_\beta$ must equal 0. Next, consider the case where $\beta_0\neq 0$. Let $t_1, t_2>0$. Then $e^{\beta_0'x}[r(t_1)-r(t_2)]$ does not depend on $x$. Since the covariance matrix of $X$ is positive definite, we can find two distinct values $x_1$ and $x_2$ of $X$ such that $e^{\beta_0'x_1}[r(t_1)-r(t_2)]=e^{\beta_0'x_2}[r(t_1)-r(t_2)]$. This implies that $r(t_1)=r(t_2)$, from which we deduce that $h_{\Lambda_j}(t)$ has to be constant (say, equal to $\alpha$) for almost every $t\in(0,\tau]$. From (\ref{inj4}), we then deduce that $h_{\Lambda_j}(0)=\alpha$, which further implies that $h_\beta=0$, $\alpha=0$, and thus $h_{\Lambda_j}(t)=0$ for almost every $t\in[0,\tau]$ ($j\in\mathcal K$). This, together with (\ref{inj2}) implies that $h_{\gamma_j}=0$, $j\in\mathcal K$.
\\
Let $k=K$. Then $\sigma_{\Lambda_K}(h)(u)=P_{\theta_0}\left[1\{S=K\}\phi(u,\mathcal O,K,\theta_0)\right]h_{\Lambda_K}(u)=0$ for all $u\in[0,\tau]$ since $h_\beta=0$ and $h_\gamma=0$. By assumptions (c) and (d), for every $u\in[0,\tau]$ and $k\in\mathcal K$,
\begin{eqnarray*}
P_{\theta_0}\left[1\{S=k\}\phi(u,\mathcal O,k,\theta_0)\right]&=&P_{\theta_0}\left[1\{S=k\}Y(u)Q(\mathcal O,k,\theta_0)e^{\beta_0'X}\right]\\
&\geq&P_{\theta_0}\left[1\{S=k\}Y(\tau)Re^{\beta_0'X}\right]>0,
\end{eqnarray*}
hence we conclude that $h_{\Lambda_K}$ is identically equal to 0 on $[0,\tau]$. Next, considering $\sigma_{\Lambda_{K-1}}(h)(u)=0$ with $h_\beta=0$, $h_\gamma=0$ and $h_{\Lambda_K}= 0$, we conclude similarly that $h_{\Lambda_{K-1}}(u)=0$ for every $u\in[0,\tau]$. It follows that $h_{\Lambda_j}$ is identically equal to 0 on $[0,\tau]$ for every $j\in\mathcal K$. Therefore, $\sigma$ is one-to-one.
\\
\\
$\Box$
\\
\\
We now turn to the proof of Theorem \ref{asympn} itself. Similar to \cite{fang05}, we get that
\begin{eqnarray*}
&&\sqrt n\left(h_\beta'(\widehat\beta_n-\beta_0)+h_\gamma'(\widehat\gamma_n-\gamma_0)+\sum_{k=1}^K\int_0^\tau h_{\Lambda_k}(s)\,d(\widehat\Lambda_{k,n}-\Lambda_{k,0})(s)\right)\\
&&\hspace{2cm}=\sqrt n\left(S_n(\theta_0)(\sigma^{-1}(h))- P_{\theta_0} \left[ S_1(\theta_0) (\sigma^{-1}(h)) \right] \right)+o_p(1),
\end{eqnarray*}
where $S_n$ is given by (\ref{score}). Consider the subset $\{(h_\beta,0,0; k\in\mathcal K)|h_\beta\in\R^p\}\subset H$ and let $\widetilde h$ be an element of this subset. Setting $h=\widetilde h$ in the above equation yields
\begin{eqnarray}\label{eqbeta}
\sqrt nh_\beta'(\widehat\beta_n-\beta_0)=\sqrt n\left(S_n(\theta_0)(\sigma^{-1}(\widetilde h))- P_{\theta_0} \left[ S_1(\theta_0) (\sigma^{-1}(\widetilde h)) \right] \right)+o_p(1).
\end{eqnarray}
By Lemma \ref{espsc0}, the central limit theorem, and Slutsky's theorem, $\sqrt nh_\beta'(\widehat\beta_n-\beta_0)$ is asymptotically normal with mean 0 and variance $P_{\theta_0}[ S_1(\theta_0) (\sigma^{-1}(\widetilde h))^2]$. If $h\in H$, direct calculation yields
\begin{eqnarray*}
&&S_1(\theta_0)(h)^2=h_\beta'S_\beta(\theta_0)S_\beta(\theta_0)'h_\beta+h_\gamma'S_\gamma(\theta_0)S_\gamma(\theta_0)'h_\gamma+2h_\beta'S_\beta(\theta_0)S_\gamma(\theta_0)'h_\gamma\\
&&\hspace{2.5cm}+2h_\beta'S_\beta(\theta_0)\left(\sum_{k=1}^K S_{\Lambda_k}(\theta_0)(h_{\Lambda_k})\right)+2h_\gamma'S_\gamma(\theta_0)\left(\sum_{k=1}^K S_{\Lambda_k}(\theta_0)(h_{\Lambda_k})\right)\\
&&\hspace{2.5cm}+\sum_{k=1}^K\left(Q(\mathcal O,k,\theta_0)\left[h_{\Lambda_k}(T)\Delta-\exp(\beta_0'X)\int_0^{T}h_{\Lambda_k}(s)\,d\Lambda_{k,0}(s)\right]\right)^2\\
&&\hspace{2.5cm}+2\sum_{k=1}^K\sum_{j>k}\left(Q(\mathcal O,k,\theta_0)\left[h_{\Lambda_k}(T)\Delta-\exp(\beta_0'X)\int_0^{T}h_{\Lambda_k}(s)\,d\Lambda_{k,0}(s)\right]\right)
\\
&&\hspace{4cm}\times
\left(Q(\mathcal O,j,\theta_0)\left[h_{\Lambda_j}(T)\Delta-\exp(\beta_0'X)\int_0^{T}h_{\Lambda_j}(s)\,d\Lambda_{j,0}(s)\right]\right).
\end{eqnarray*}
Taking expectation followed by some tedious algebraic manipulations and re-arrangement of terms yield that
\begin{eqnarray*}
P_{\theta_0} \left[S_1(\theta_0)(h)^2 \right]=h_\beta'\sigma_\beta(h)+h_\gamma'\sigma_\gamma(h)+\sum_{k=1}^K\int_0^\tau\sigma_{\Lambda_k}(h)(u)h_{\Lambda_k}(u)\,d\Lambda_{k,0}(u).
\end{eqnarray*}
Therefore
\begin{eqnarray*}
P_{\theta_0} \left[ S_1(\theta_0) (\sigma^{-1}(\widetilde h))^2\right]&=&\sigma_\beta^{-1}(\widetilde h)'\sigma_\beta(\sigma^{-1}(\widetilde h))+\sigma_\gamma^{-1}(\widetilde h)'\sigma_\gamma(\sigma^{-1}(\widetilde h))\\
&&\hspace{2cm}+\sum_{k=1}^K\int_0^\tau \sigma_{\Lambda_k}^{-1}(\widetilde h)(u) \sigma_{\Lambda_k} (\sigma^{-1}(\widetilde h))(u)\,d\Lambda_{k,0}(u)\\
&=&h_\beta'\sigma_\beta^{-1}(\widetilde h),
\end{eqnarray*}
where the last equality comes from the fact that
\begin{eqnarray*}
\sigma(\sigma^{-1}(\widetilde h))=(\sigma_\beta (\sigma^{-1}( \widetilde h)), \sigma_\gamma (\sigma^{-1}(\widetilde h)), \sigma_{\Lambda_k}(\sigma^{-1}(\widetilde h));k\in\mathcal K)=\widetilde h.
\end{eqnarray*}
Now, recall that the linear map $\widetilde \sigma_\beta^{-1}:\R^p\rightarrow \R^p$ was defined in Section \ref{asymprp} as a restricted version of $\sigma_\beta^{-1}$, by setting $\widetilde \sigma_\beta^{-1} (h_\beta)= \sigma_\beta^{-1} (\widetilde h)$ for any $\widetilde h$ of the form $(h_\beta, 0, 0; k\in\mathcal K)$. Let $\{e_1,\ldots,e_p\}$ be the canonical basis of $\R^p$ and $\Sigma_\beta=(\widetilde \sigma_\beta^{-1}(e_1),\ldots,\widetilde \sigma_\beta^{-1}(e_p))$. Then for any $h_\beta\in\R^p$, we have $\widetilde \sigma_\beta^{-1}(h_\beta)=\Sigma_\beta h_\beta$ and thus $P_{\theta_0}[ S_1(\theta_0) (\sigma^{-1}(\widetilde h))^2]=h_\beta'\Sigma_\beta h_\beta$. Hence, for every $h_\beta\in\R^p$, $\sqrt n h_\beta'(\widehat\beta_n-\beta_0)$ converges in distribution to $\mathcal N(0,h_\beta'\Sigma_\beta h_\beta)$. By the Cram\'er-Wold device, $\sqrt n(\widehat\beta_n-\beta_0)$ converges in distribution to $\mathcal N(0,\Sigma_\beta)$.
\\\\
Now, for $j=1\ldots,p$, denote $\widetilde h_j=(e_j,0,0; k\in\mathcal K)$. Letting $h=\widetilde h_j$ for each $j=1\ldots,p$ in turn in (\ref{eqbeta}) yields
\begin{eqnarray*}
\sqrt n(\widehat\beta_n-\beta_0)=\frac{1}{\sqrt n}\sum_{i=1}^n l_\beta(\mathcal O_i, \theta_0)+o_p(1),
\end{eqnarray*}
where
\begin{eqnarray*}
l_\beta(\mathcal O, \theta_0)=\Sigma_\beta S_\beta(\theta_0)+\Xi S_\gamma(\theta_0)+\sum_{k=1}^K S_{\Lambda_k}(\theta_0)(\Xi^\ast),
\end{eqnarray*}
$\Xi$ and $\Xi^\ast$ are $(p\times q)$ and $(p\times 1)$ matrices respectively defined by
\begin{eqnarray*}
\Xi=\left(
\begin{array}{c}
\sigma_\gamma^{-1}(\widetilde h_1)'\\
\vdots\\
\sigma_\gamma^{-1}(\widetilde h_p)'
\end{array}\right) \quad\mbox{ and }\quad 
\Xi^\ast=\left(
\begin{array}{c}
\sigma_{\Lambda_k}^{-1}(\widetilde h_1)\\
\vdots\\
\sigma_{\Lambda_k}^{-1}(\widetilde h_p)
\end{array}\right),
\end{eqnarray*}
and $S_{\Lambda_k}(\theta_0)$ is applied componentwise to $\Xi^\ast$. Thus $\widehat\beta_n$ is an asymptotically linear estimator for $\beta_0$, and its influence function $l_\beta(\mathcal O, \theta_0)$ belongs to the tangent space spanned by the score functions. It follows that $l_\beta(\mathcal O, \theta_0)$ is the efficient influence function for $\beta_0$, and that $\widehat\beta_n$ is semiparametrically efficient (see \cite{bkrw} or \cite{tsi06}).
\\
\\
$\Box$
\\
\\
\textbf{A.4 Proof of Theorem \ref{asympnpn}}
\\
\\
The proof of asymptotic normality of $\sqrt n(\widehat\gamma_n-\gamma_0)$ proceeds along the same line as for $\sqrt n(\widehat\beta_n-\beta_0)$, and is therefore omitted.
\\
Next, for any $t\in(0,\tau)$ and $j\in\mathcal K$, the asymptotic normality of $\sqrt n(\widehat\Lambda_{j,n}(t)-\Lambda_{j,0}(t))$ can be proved by using a similar argument with $\widetilde h$ replaced by $h_{(j,t)}=(h_\beta, h_\gamma, h_{\Lambda_k}; k\in\mathcal K)$, where $h_\beta=0$, $h_\gamma=0$, $h_{\Lambda_j}(\cdot)=1\{\cdot\leq t\}$ ($t\in(0,\tau)$ and $j\in\mathcal K$), and $h_{\Lambda_k}=0$ for every $k\in\mathcal K, k\neq j$. Details are omitted.
\\
\\
$\Box$
\\
\\
\textbf{A.5 Proof of Theorem \ref{varasymp}}
\\\\
The proof of Theorem \ref{varasymp} parallels the proof of Theorem 3 in \cite{parner} and thus, will be kept brief. Let $\widehat \sigma_n=(\widehat \sigma_{\beta,n}, \widehat \sigma_{\gamma,n}, \widehat \sigma_{\Lambda_k,n}; k \in \mathcal K)$ be defined as  $\sigma$ with all of the $\theta_0$ and $P_{\theta_0}$ replaced by $\widehat\theta_n$ and $\mathbb P_n$ respectively. Similar to the proof of Theorem 3 in \cite{parner}, it can be shown that $\widehat\sigma_n$ converges in probability to $\sigma$ uniformly over $H$ and that its inverse $\widehat \sigma_n^{-1}=(\widehat \sigma_{\beta,n}^{-1}, \widehat \sigma_{\gamma,n}^{-1}, \widehat \sigma_{\Lambda_k,n}^{-1}; k \in \mathcal K)$ is such that $\widehat \sigma_n^{-1}(h)$ converges to $\sigma^{-1}(h)$ in probability.

For every $h_\beta$, the asymptotic variance of $\sqrt n h_\beta'(\widehat\beta_n-\beta_0)$ is $h_\beta'\sigma_\beta^{-1} ((h_\beta,0, 0; k\in \mathcal K))$, which is consistently estimated by $h_\beta'\widehat \sigma_{\beta,n}^{-1} ((h_\beta, 0, 0; k\in \mathcal K))$. Let $h_n=(h_{\beta,n}, h_{\gamma,n}, h_{\Lambda_k,n}; k\in \mathcal K)=\widehat \sigma_n^{-1}((h_\beta,0, 0; k\in \mathcal K))$. Then $\widehat \sigma_n(h_n)=(h_\beta,0, 0; k\in \mathcal K)$, or
\begin{eqnarray}\label{sysbeta}
\left\{
\begin{array}{l}
\widehat \sigma_{\beta,n}(h_n)=h_\beta\\
\widehat \sigma_{\gamma,n}(h_n)=0\\
\widehat \sigma_{\Lambda_k,n}(h_n)(u)=0, \quad k\in\mathcal K, \quad u\in[0,\tau].\\
\end{array}\right.
\end{eqnarray}
In particular, letting $u=T_1,\ldots,T_n$ in (\ref{sysbeta}) yields the following system of equations:
\begin{eqnarray}\label{sys}
\mathbb A_n\left(
\begin{array}{c}
h_{\beta,n}\\
h_{\gamma,n}\\
h_{\Lambda,n}
\end{array}
\right)=\left(
\begin{array}{l}
h_\beta\\
0_q\\
0_{Kn}
\end{array}
\right),
\end{eqnarray}
where $h_{\Lambda,n}=(h_{\Lambda_1,n}(T_1), \ldots,h_{\Lambda_1,n}(T_n),\ldots, h_{\Lambda_K,n}(T_1), \ldots, h_{\Lambda_K,n}(T_n))'$, and $\mathbb A_n$ is defined by (\ref{bm}). Some simple algebra on (\ref{sys}) yields that $h_{\beta,n}=\widehat\Sigma_{\beta,n} h_\beta$ where $\widehat\Sigma_{\beta,n}$ is defined in Section \ref{asymprp}, and therefore $h_\beta' \widehat\Sigma_{\beta,n} h_\beta$ is a consistent estimator of the asymptotic variance of $\sqrt n h_\beta'(\widehat\beta_n-\beta_0)$ for every $h_\beta$. We conclude that $\widehat\Sigma_{\beta,n}$ converges in probability to $\Sigma_\beta$. The consistency of $\widehat\Sigma_{\gamma,n}$ proceeds along the same lines and is therefore omitted.

We now turn to the estimation of the asymptotic variance of $\widehat \Lambda_{j,n}(t)$, for $t\in (0,\tau)$ and $j\in\mathcal K$. By the dominated convergence theorem and the consistency of $\widehat \sigma_n^{-1}$,
\begin{eqnarray*}
\int_0^t\widehat \sigma_{\Lambda_j,n}^{-1}(h_{(j,t)})(u)\,d\widehat\Lambda_{j,n}(u)
\end{eqnarray*}
converges to $v_j^2(t)=\int_0^t \sigma_{\Lambda_j}^{-1}(h_{(j,t)})(u)\, d\Lambda_{j,0}(u)$, where we recall that $h_{(j,t)}$ is the element $(h_\beta, h_\gamma, h_{\Lambda_k}; k\in\mathcal K)$ such that $h_\beta=0$, $h_\gamma=0$, $h_{\Lambda_j}(\cdot)=1\{\cdot\leq t\}$ for some $t\in(0,\tau)$ and $j\in\mathcal K$, and $h_{\Lambda_k}=0$ for every $k\in\mathcal K, k\neq j$. Letting $\widetilde h_n=(\widetilde h_{\beta,n}, \widetilde h_{\gamma,n}, \widetilde h_{\Lambda_k,n}; k\in \mathcal K)=\widehat \sigma_n^{-1}(h_{(j,t)})$, we get that $\widehat \sigma_n(\widetilde h_n)=h_{(j,t)}$ or:
\begin{eqnarray}\label{syslambda}
\left\{
\begin{array}{l}
\widehat \sigma_{\beta,n}(\widetilde h_n)=0\\
\widehat \sigma_{\gamma,n}(\widetilde h_n)=0\\
\widehat \sigma_{\Lambda_j,n}(\widetilde h_n)(u)=1\{u\leq t\}, \quad u\in[0,\tau]\\
\widehat \sigma_{\Lambda_k,n}(\widetilde h_n)(u)=0, \quad k\in\mathcal K, \quad k\neq j, \quad u\in[0,\tau].\\
\end{array}\right.
\end{eqnarray}
In particular, letting $u=T_1,\ldots,T_n$ in (\ref{syslambda}) yields the system of equations
\begin{eqnarray*}
\mathbb A_n\left(
\begin{array}{c}
\widetilde h_{\beta,n}\\
\widetilde h_{\gamma,n}\\
\widetilde h_{\Lambda,n}
\end{array}
\right)=\left(
\begin{array}{l}
0_p\\
0_q\\
U_{(j,t)}^n
\end{array}
\right),
\end{eqnarray*}
with the notations $\widetilde h_{\Lambda,n}=(\widetilde h_{\Lambda_1,n}(T_1), \ldots,\widetilde h_{\Lambda_1,n}(T_n),\ldots, \widetilde h_{\Lambda_K,n}(T_1), \ldots, \widetilde h_{\Lambda_K,n}(T_n))'$ and $U_{(j,t)}^n=(0_{(j-1)n}',1\{T_1\leq t\},\ldots,1\{T_n\leq t\},0_{(K-j)n}')'$. Similar algebra as above yields
\begin{eqnarray*}
\widetilde h_{\Lambda,n}=\widehat\Sigma_{\Lambda,n} U_{(j,t)}^n,
\end{eqnarray*}
where $\widehat\Sigma_{\Lambda,n} $ is defined in Section \ref{asymprp}. Now, $\int_0^t\widehat \sigma_{\Lambda_j,n}^{-1}(h_{(j,t)})(u)\,d\widehat\Lambda_{j,n}(u)$ verifies
\begin{eqnarray*}
\int_0^t\widehat \sigma_{\Lambda_j,n}^{-1}(h_{(j,t)})(u)\,d\widehat\Lambda_{j,n}(u)&=&\sum_{i=1}^n \widehat \sigma_{\Lambda_j, n}^{-1}(h_{(j,t)})(T_i)\widehat{\Delta \Lambda_{j,n}}(T_i) 1\{T_i\leq t\} \\
&=& \widehat\Xi_{(j,t)}^{n'}\widetilde h_{\Lambda,n},
\end{eqnarray*}
where $\widehat\Xi_{(j,t)}^n=\left(0_{(j-1)n}',\widehat{\Delta \Lambda_{j,n}}(T_1) 1\{T_1\leq t\},\ldots,\widehat{\Delta \Lambda_{j,n}}(T_n) 1\{T_n\leq t\},0_{(K-j)n}'\right)'$. It follows that $\widehat\Xi_{(j,t)}^{n'}\widehat\Sigma_{\Lambda,n} U_{(j,t)}^n$ is a consistent estimator of $v_j^2(t)$, which concludes the proof.
\\
\\
$\Box$
\bibliographystyle{Chicago}
\bibliography{biblio}
\end{document}